\documentclass[12pt]{article}
\usepackage{amsthm,amsmath,amssymb,amscd,verbatim}
\usepackage{graphics,pifont}
\usepackage{amssymb}
\usepackage{graphicx}
\date{}
\newcommand{\f}{\frac}
\newcommand{\ov}{\overline}
\newcommand{\ra}{\rightarrow}

\begin{document}

\title{{\bf An example of unbounded chaos}}
\author{Bau-Sen Du  \\
Institute of Mathematics \\ Academia Sinica \\ Taipei 10617, Taiwan \\
dubs@math.sinica.edu.tw\\}
\maketitle

\begin{abstract} 
\noindent
Let $\phi(x) = |1 - \frac 1x|$ for all $x > 0$.  Then we extend $\phi(x)$ in the usual way to become a continuous map from the compact topological (but not metric) space $[0, \infty]$ onto itself which also maps the set of irrational points in $(0, \infty)$ onto itself.  In this note, we show that (1) on $[0, \infty]$, $\phi(x)$ is topologically mixing, has dense irrational periodic points, and has topological entropy $\log \lambda$, where $\lambda$ is the unique positive zero of the polynomial $x^3 - 2x -1$; (2) $\phi(x)$ has bounded uncountable {\it invariant} 2-scrambled sets of irrational points in $(0, 3)$; (3) for any countably infinite set $X$ of points (rational or irrational) in $(0, \infty)$, there exists a dense unbounded uncountable {\it invariant} $\infty$-scrambled set $Y$ of irrational transitive points in $(0, \infty)$ such that, for any $x \in X$ and any $y \in Y$, we have $\limsup_{n \to \infty} |\phi^n(x) - \phi^n(y)| = \infty$ and $\liminf_{n \to \infty} |\phi^n(x) - \phi^n(y)| = 0$.  This demonstrates the true nature of chaos for $\phi(x)$ (see \cite{du2, v1, v2}).
\medskip


\end{abstract}

\noindent
{\bf 1. Introduction}

\noindent
Let $I$ denote an interval in the real line and let $f$ be a continuous map from $I$ into itself.  It is well-known that if $f$ has a periodic point of period not an integral power of 2 then there exist a positive number $\delta$ and an uncountable set $S$ such that $$\text{for any} \,\, x \ne y \,\, \text{in} \,\, S, \,\, \limsup_{n \to \infty} |f^n(x) - f^n(y)| \ge \delta \,\,\, \text{and} \,\,\, \liminf_{n \to \infty} |f^n(x) - f^n(y)| = 0. \qquad (*)$$ Such a set $S$ is called a $\delta$-scrambled set of $f$.  When the inequality in (*) holds for all positive numbers $\delta$ (we follow the convention that $|\pm \infty \pm$ any real number$| = \infty$ and $\infty >$ any real number), we call such set $S$ an $\infty$-scrambled set.  Can this scrambled set $S$ be chosen to be invariant under $f$?  That is, can $f(S) \subset S$?  The answer is no in general.  This is because if $c$ is a point in $I$ such that $\liminf_{n \to \infty} |f^n(c) - f^n(f(c))| = 0$, then the $\omega$-limit set $\omega(f, c)$ of $c$ must contain a fixed point of $f$, where the $\omega$-limit set $\omega(f, c)$ of a point $c$ is the set of all points $x$ with the property that there is an increasing sequence $<n_i>$ of positive integers such that $\lim_{i \to \infty} f^{n_i}(c) = x$.  This fact is useful in constructing examples of maps without {\it invariant} scrambled sets.  For example, let $g$ be the continuous map from $[0, 1]$ onto itself such that (i) $g(0)= 1$, $g(1) = 1/2$, and $g(1/2) = 0$; (ii) $g(1/6) = 1/3$ and $g(1/3) = 1/6$; and (iii) $g$ is linear on each of the intervals $[0, 1/6]$, $[1/6, 1/3]$, $[1/3, 1/2]$ and $[1/2, 1]$.  Then the point $x = 1/4$ is the unique fixed point of $g$ and every point in $(1/6, 1/4) \cup (1/4, 1/3)$ is a period-2 point of $g$.  Therefore, $g$ has no {\it invariant} scrambled sets although $g$ has the period-3 orbit $\{0, 1/2, 1 \}$.  On the other hand, if $f$ is turbulent, i.e., if there exist two compact subintervals $I_0$ and $I_1$ of $I$ with at most one point in common such that $f(I_0) \cap f(I_1) \supset I_0 \cup I_1$, then we can find such an {\it invariant} scrambled set \cite{du1}.  

\begin{figure}[ht] 
\begin{center}
\includegraphics[width=10cm,height=5cm]{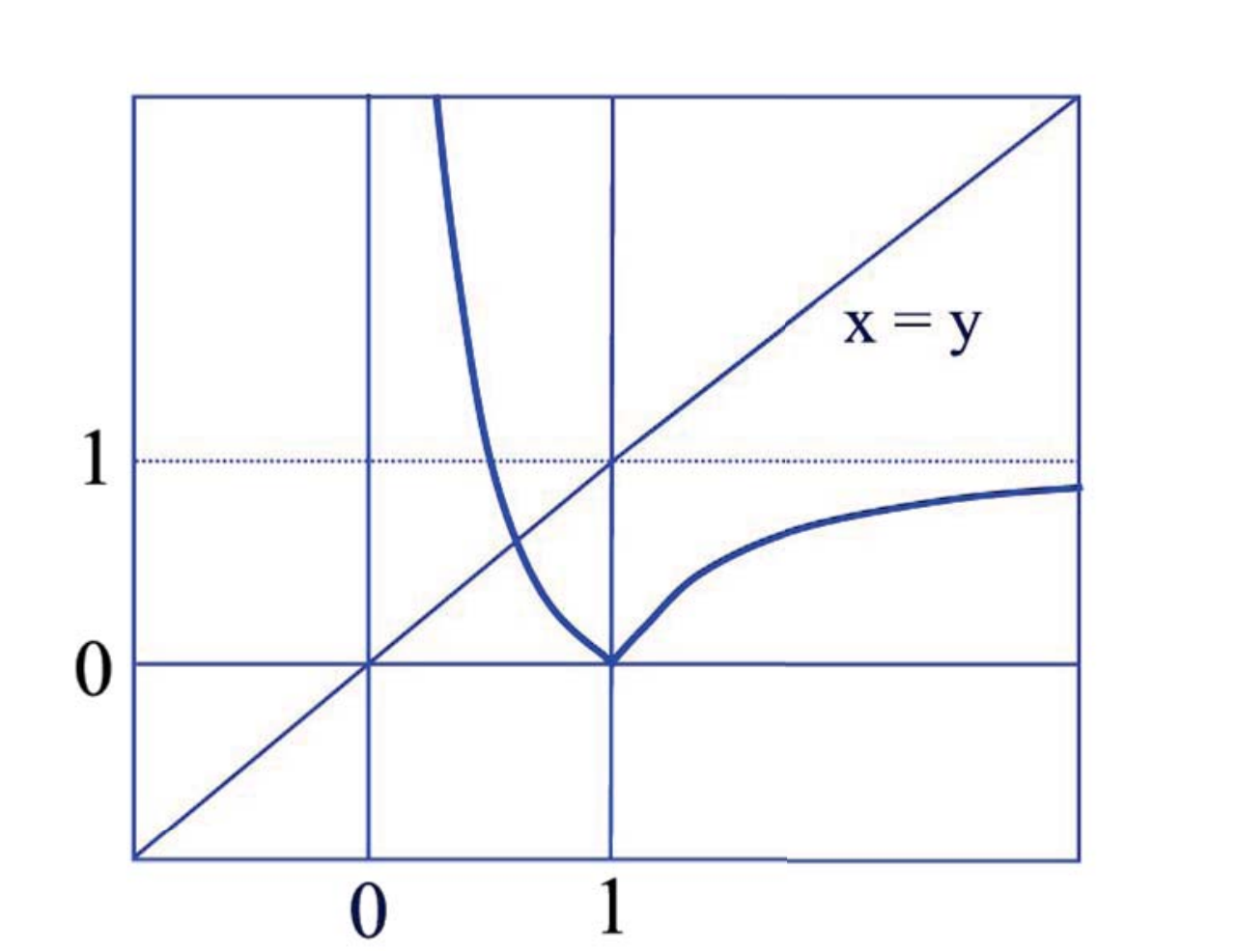} 
\caption{The graph of the map $\phi(x) = |1 - {\frac 1x}|, x > 0$.}
\end{center}
\end{figure}

In studying the periodicity and convergence of the second-order difference equation $x_{n+1}$ $= |x_n - x_{n-1}|$, Sedaghat \cite{se1, se2} introduced the following interval map, see Figure 1, $$\phi(x) = \biggl|1 - \frac 1x \biggr| \,\,\, \text{defined on} \,\,\, (0, \infty)$$ and showed, among other things, that $\phi(x)$ has  uncountable scrambled sets.  Following the work of Sedaghat \cite{se1, se2}, we may ask the following two questions:
\begin{itemize}
\item[(1)] Does $\phi(x)$ have bounded uncountable {\it invariant} scrambled sets ?

\item[(2)] Does $\phi(x)$ have unbounded uncountable {\it invariant} scrambled sets ?
\end{itemize}

In this note, we shall answer both questions afirmatively.  We shall also show that for any countably infinite set $X$ of points (rational or irrational) in $(0, \infty)$, there exists a dense unbounded uncountable invariant $\infty$-scrambled set $Y$ of irrational transitive points (transitive points are points with dense orbits) in $(0, \infty)$ such that for any $x \in X$ and any $y \in Y$, we have $\limsup_{n \to \infty} |\phi^n(x) - \phi^n(y)| = \infty$ and $\liminf_{n \to \infty} |\phi^n(x) - \phi^n(y)| = 0$.  For these purposes, we shall use symbolic dynamics.
\bigskip

\noindent
{\bf 2. Symbolic dynamics}

\noindent
Let $\Sigma_2 = \{ \beta : \beta = \beta_0\beta_1 \cdots$, where $\beta_i = 0 \, \text{or} \, 1 \}$ be the metric space with metric $d$ defined by $$d(\beta_0\beta_1 \cdots, \gamma_0\gamma_1 \cdots) = \sum_{i=0}^\infty \frac {|\beta_i-\gamma_i|}{2^{i+1}}.$$ Let $\sigma$ be the shift map defined by $\sigma(\beta_0\beta_1\beta_2 \cdots) = \beta_1\beta_2 \cdots$.  Then $\sigma$ is a continuous, two-to-one map from $\Sigma_2$ onto itself.  In the sequel, for any finite sequences $\beta_0\beta_1\beta_2\cdots \beta_k$, $k \ge 2$, of 0's and 1's, we also define $\sigma(\beta_0\beta_1\beta_2\cdots \beta_k) = \beta_1\beta_2\cdots \beta_k$.  

\begin{figure}[ht] 
\begin{center}
\includegraphics[width=16cm,height=8cm]{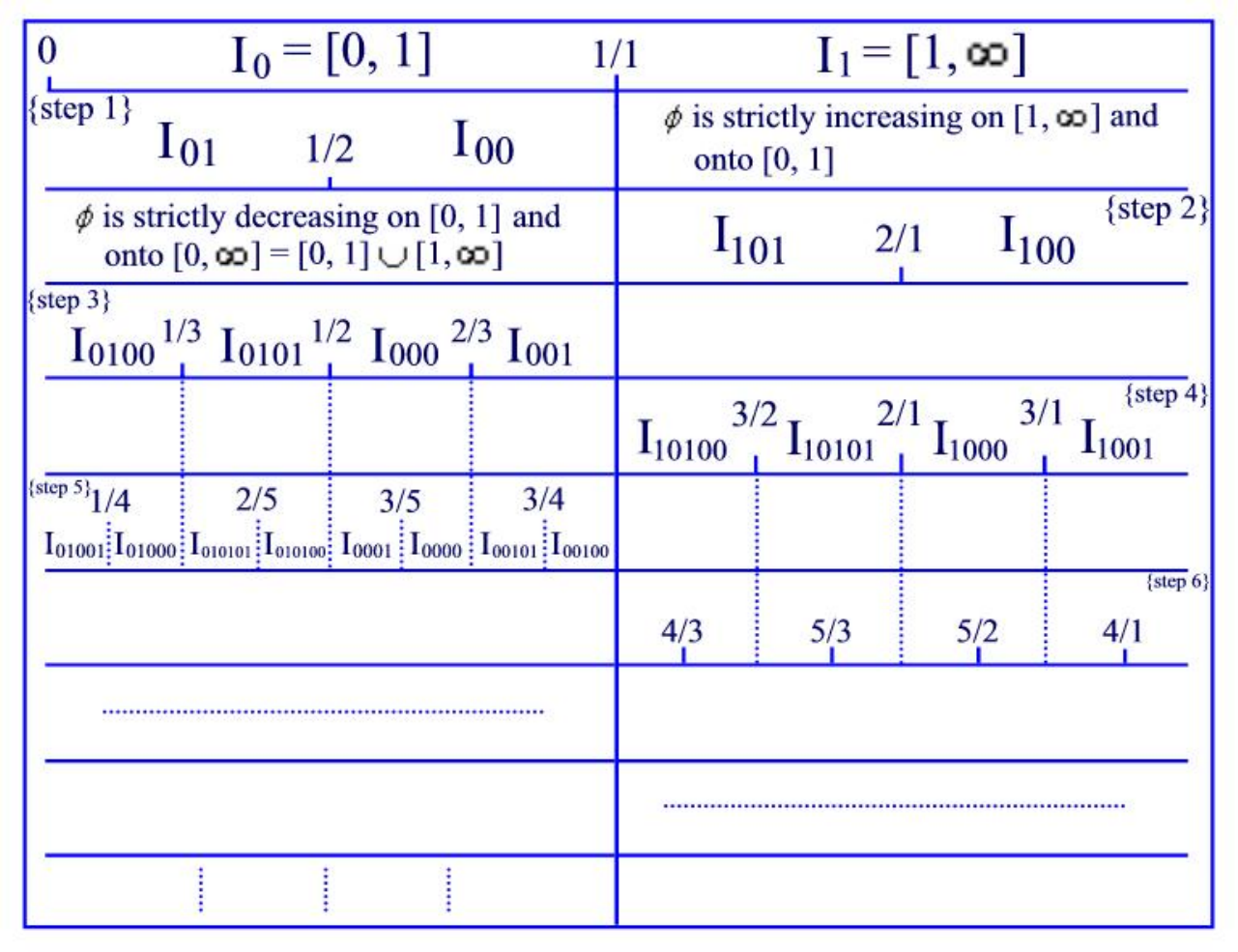} 
\caption{Recursive procedures to obtain $I_{\gamma_0\gamma_1\gamma_2 \cdots \gamma_{n_i-1}\gamma_{n_i}}$'s.}
\end{center}
\end{figure}

Let $\phi(x) = |1 - 1/x|$ for all $x > 0$.  We now compactify the space $[0, \infty)$ by including the symbol $\infty$ and define $\phi(0) = \infty$ and $\phi(\infty) = 1$.  Then $\phi$ is a continuous map from the compact topological (but not metric) space $[0, \infty]$ onto itself.  Let $I_0 = [0, 1]$ and $I_1 = [1, \infty]$.  Then $\phi(I_0) = [0, \infty] = I_0 \cup I_1$ and $\phi(I_1) = [0, 1] = I_0$.  Since $\phi(I_0) = I_0 \cup I_1$ and $\phi$ is strictly decreasing on $I_0$, we can split $I_0$ up into two compact subintervals $I_{00}$ and $I_{01}$ with one point in common such that $\phi(I_{0i}) = I_i = I_{\sigma(0i)}$, $i = 0, 1$.  Since $\phi(I_1) = I_0$ which is disjoint from the interior of $I_1$, we cannot have two proper compact subintervals $I_{10}$ and $I_{11}$ of $I_1$ such that $\phi(I_{1i}) = I_i$, $i = 0, 1$.  However, since $\phi(I_1) = I_0 = I_{00} \cup I_{01}$ and $\phi$ is strictly increasing on $I_1$, we can split $I_1$ up into two compact subintervals $I_{100}$ and $I_{101}$ with one point in common such that $\phi(I_{10i}) = I_{0i} = I_{\sigma(10i)}$, $i = 0, 1$.  Therefore, although we {\it cannot} split $I_1$ up into two compact subintervals $I_{10}$ and $I_{11}$ with one point in common such that $\phi(I_{1i}) = I_i, i = 0, 1$, we can split $I_1$ up into two compact subintervals $I_{100}$ and $I_{101}$ (subscripts obtained by inserting the number 1 before each of $00$ and $01$) with one point in common such that $\phi(I_{10i}) = I_{0i}, i = 0, 1$.  Similarly, we can continue with the above procedures (see Figure 2) as follows: 
\begin{itemize}
\item[(1)] since $\phi$ is strictly decreasing on $I_0 = [0, 1]$ and $\phi(I_0) = [0, \infty] = I_{100} \cup I_{101} \cup I_{00} \cup I_{01}$ (in that order which is obtained from Steps 2 and 1 or, $(2n)^{th}$ and $(2n-1)^{st}$ steps in Figure 2 from right to left), we can split $I_0$ up, from left to right, into four compact subintervals with pairwise disjoint interiors, $I_{0100}, I_{0101}, I_{000}, I_{001}$ (obtained by inserting a 0 before the "previous" four subscripts) such that $I_{\gamma_0\gamma_1 \cdots \gamma_{n_i}} \subset I_{\gamma_0\gamma_1 \cdots \gamma_{n_i-1}}$ and $\phi(I_{\gamma_0\gamma_1 \cdots \gamma_{n_i}}) = I_{\gamma_1\gamma_2 \cdots \gamma_{n_i}}$ wherever defined, and the union of all these intervals is $[0, 1]$,    

\item[(2)] since $\phi$ is strictly increasing on $I_1 = [1, \infty]$ and $\phi(I_1) = [0, 1] = I_{0100} \cup I_{0101} \cup I_{000} \cup I_{001}$ (in that order which is obtained from Procedure 1, i.e., Step 3 or $(2n+1)^{st}$ step in Figure 2 from left to right), we can split $I_1$ up, from left to right, into four compact subintervals with pairwise disjoint interiors, $I_{10100}, I_{10101}, I_{1000}, I_{1001}$ (obtained by inserting a 1 before the "previous" four subscripts) such that $I_{\gamma_0\gamma_1 \cdots \gamma_{n_i}} \subset I_{\gamma_0\gamma_1 \cdots \gamma_{n_i-1}}$ and $\phi(I_{\gamma_0\gamma_1 \cdots \gamma_{n_i}}) = I_{\gamma_1\gamma_2 \cdots \gamma_{n_i}}$ wherever defined, and the union of all these intervals is $[1, \infty]$,   

\item[(3)] repeat Procedure 1 first and then Procedure 2 indefinitely.  
\end{itemize}

\noindent
Consequently, $I_0 = [0, 1], I_1 = [1, \infty]$, $I_{00} = [1/2, 1]$, $I_{01} = [0, 1/2]$ and, for any finite sequence $\gamma_0\gamma_1\gamma_2 \cdots \gamma_{n_i} = \gamma_0\gamma_1\gamma_2 \cdots \gamma_{n_i-2}00$ or $\gamma_0\gamma_1\gamma_2 \cdots \gamma_{n_i-2}01$, $n_i \ge 2$, of 0's and 1's such that if $\gamma_j = 1$ for some $0 \le j \le n_i-2$ then $\gamma_{j+1} = 0$, the compact interval $I_{\gamma_0\gamma_1\gamma_2 \cdots \gamma_{n_i}}$ is defined so that $I_{\gamma_0\gamma_1\gamma_2\cdots \gamma_{n_i}} \subset I_{\gamma_0\gamma_1\gamma_2\cdots \gamma_{n_i-1}}$ and $\phi(I_{\gamma_0\gamma_1\gamma_2\cdots \gamma_{n_i}}) = I_{\gamma_1\gamma_2\cdots \gamma_{n_i}}$.  Note that, in Figure 2, if $I_{0\gamma_1\gamma_2 \cdots \gamma_{n_i}}$ is defined at $(2n-1)^{st}$ step, then $I_{10\gamma_1\gamma_2 \cdots \gamma_{n_i}}$ is defined at $(2n)^{th}$ step, and $I_{00\gamma_1\gamma_2 \cdots \gamma_{n_i}}$ is defined at $(2n+2)^{nd}$ step, while if $I_{1\gamma_1\gamma_2 \cdots \gamma_{n_i}}$ is defined at $(2n)^{th}$ step (and so $\gamma_1 = 0$), then $I_{01\gamma_1\gamma_2 \cdots \gamma_{n_i}}$ is defined at $(2n+1)^{st}$ step.  Furthermore, at the $(2n-1)^{st}$ \,($(2n)^{th}$ respectively) step, the interval $[0, 1]$ \, ($[1, \infty]$ respectively) is split up into the union of $2^n$ compact subintervals with pairwise disjoint interiors and each of these $2^n$ intervals is split up at the next step into the union of two compact subintervals with disjoint interiors.  Let $\Gamma = \{ \gamma = \gamma_0\gamma_1\gamma_2 \cdots \in \Sigma_2 : $ if $\gamma_j = 1$ then $\gamma_{j+1} = 0 \}$.  Then $\Gamma$ is a subshift of finite type and, for any $\gamma = \gamma_0\gamma_1\gamma_2 \cdots \ne \overline {100}$ in $\Gamma$ ($\ov {100}$ stands for the sequence of repeating $100$ in $\Gamma$ and note that $I_{\ov {100}} = \{ \infty \}$), the sequence $< I_{\gamma_0\gamma_1\gamma_2 \cdots \gamma_{n_i}} >_{i \ge 1}$, wherever $I_{\gamma_0\gamma_1\gamma_2 \cdots \gamma_{n_i}}$ is defined, is a nested sequence of compact intervals in $[0, \infty)$.  Therefore, the set $I_\gamma = \cap_{i \ge 0} I_{\gamma_0\gamma_1\gamma_2 \cdots \gamma_{n_i}}$ is either a nontrivial compact interval or consists of exactly one point and $\phi(I_\gamma) = I_{\sigma\gamma}$ for each $\gamma \in \Gamma$.  However, since each rational point in $[0, \infty)$ is mapped to the point 0 after a finite number of iterations of $\phi$ \cite{se1}, it is clear that each $I_\gamma$ actually consists of exactly one point.  {\it In the sequel}, we write this point as $x_\gamma$.  So, $I_\gamma = \{ x_\gamma \}$.  Since $\phi(I_\gamma) = I_{\sigma\gamma}$ for each $\gamma \in \Gamma$, we obtain that if $\sigma^m(\gamma) = \gamma$ for some $\gamma \in \Gamma$, then $\phi^m(x_\gamma) = x_\gamma$.  We note in passing that this correspondence $\gamma \to x_\gamma$ from $\Gamma$ to $[0, \infty]$ is not one-to-one, but is {\it onto}.  For each rational number $s > 0$, there are exactly two distinct $\gamma^{(1)}, \gamma^{(2)}$ in $\Gamma$ such that $I_{\gamma^{(1)}} = I_{\gamma^{(2)}} = \{ s \}$.  For example, $I_{0\overline {010}} = \{ 1 \} = I_{1 \overline {010}}$.  However, we don't need this information later on.  

We remark that, since $\phi(x)$ is a very special map, we can actually compute these compact intervals $I_{\gamma_0\gamma_1 \cdots \gamma_{n_i}}$ explicitly.  They are related to Farey sequences in the way that, if $I_0 = [0, 1] = [{\frac 01}, {\frac 11}], I_1 = [1, \infty] = [{\frac 11}, {\frac 10}]$, and if $a \ge 1, b \ge 0, c \ge 0, d \ge 1$ are integers such that $[{\frac ba}, {\frac dc}]$ (in lowest terms) is one of these compact intervals $I_{\gamma_0\gamma_1 \cdots \gamma_{n_i}}$ in $[0, \infty]$ obtained at the $n^{th}$ step, then $[{\frac ba}, {\frac {b+d}{a+c}}]$ and $[{\frac {b+d}{a+c}}, {\frac dc}]$ are two adjacent compact intervals in $[0, \infty]$ obtained at $(n+1)^{st}$ step (see Figure 2).  However, we don't need this information to achieve our goal.  
\bigskip

\noindent
{\bf 3. The extended map $\phi(x)$ is topologically mixing, has dense transitive irrational points and has dense irrational periodic points in $[0, \infty]$}

\noindent
Let $V$ be a nonempty open interval in $(0, \infty)$, then there exist an irrational point $w$ in $V$ and an element $\gamma = \gamma_0\gamma_1\gamma_2 \cdots$ in $\Gamma$ such that $I_\gamma = \{ w \}$.  Since there exists a strictly increasing sequence $<n_i>$ of positive integers such that $V \supset \{ w \} = I_{\gamma} = \cap_{i \ge 0} I_{\gamma_0\gamma_1\gamma_2 \cdots \gamma_{n_i}}$, there is an integer $k \ge 5$ such that $w \in I_{\gamma_0\gamma_1\gamma_2 \cdots \gamma_{n_k}} \subset V$.  Consequently, $I_0 = \phi^{n_k+1}(I_{\gamma_0\gamma_1\gamma_2\cdots \gamma_{n_k}00})  \subset \phi^{n_k+1}(V)$ and so $[0, \infty] \supset \phi^{n_k+2}(V) \supset \phi(I_0) = [0, \infty]$.  Thus, $\phi^{n_k+2}(V) = [0, \infty]$.  This implies trivially that the extended map $\phi(x)$ is topologically mixing on $[0, \infty]$.  Later in section 6, we shall see that $\phi(x)$ has dense transitive irrational points in $[0, \infty]$.  Finally, since each rational point in $(0, \infty)$ is mapped into the period-3 orbit $\{1, 0, \infty \}$ under a finite number of iterations of $\phi$ and since $\{ 1 \} = I_{\overline{001}} = I_{101\overline{001}}$, the unique point $p$ in the set $I_{\overline {\gamma_0\gamma_1\gamma_2\cdots \gamma_{n_k} \, 000}} \subset I_{\gamma_0\gamma_1\gamma_2\cdots \gamma_{n_k}000} \subset I_{\gamma_0\gamma_1\gamma_2\cdots \gamma_{n_k}} \subset V$ is an irrational periodic point of $\phi(x)$ in $V$.  This shows that the irrational periodic points of $\phi(x)$ are dense in $(0, \infty)$.    
\bigskip

\noindent
{\bf 4. The topological entropy of the extended map $\phi(x)$ on $[0, \infty]$}

\noindent
Since the extended map $\phi(x)$ is continuous on the compact topological space $[0, \infty]$, we can compute its topological entropy.  We first define the continuous piecewise linear map $f : [0, 1] \to [0, 1]$ by putting $f(x) = 1 - 2x$ for $0 \le x \le 1/2$ and $f(x) = x - 1/2$ for $1/2 \le x \le 1$.  Then $\{ 0, \frac 12, 1 \}$ is a period-3 orbit of $f$ and it is well-known that the topological entropy of $f$ is $\log \lambda$, where $\lambda$ is the unique positive zero of the polynomial $x^3 - 2x - 1$.  In this section, we shall show that $\phi$ is topologically conjugate to the map $f$ and so $\phi$ also has entropy $\log \lambda$.  For this purpose, we define modified Farey sequences on $[0, \infty]$ inductively as follows: $F_0 = \{ \frac 01, \frac 10 \}$ and if $\frac ab$ and $\frac cd$ (both are in lowest terms) are two consecutive fractions of $F_n$, then $\frac ab, \frac {a+c}{b+d}, \frac cd$ are consecutive fractions of $F_{n+1}$.  The first four modified Farey sequences are $$F_0 = \{ \frac 01, \frac 10 \}; \quad F_1 = \{ \frac 01, \frac 11, \frac 10 \}; \quad F_2 = \{ \frac 01, \frac 12, \frac 11, \frac 21, \frac 10 \}; \quad F_3 = \{ \frac 01, \frac 13, \frac 12, \frac 23, \frac 11, \frac 32, \frac 21, \frac 31, \frac 10 \}.$$  Furthermore, if we write, for each integer $n \ge 1$, $F_n = \{ 0/1 = a_{0, n}/b_{0, n}$, $\, a_{1, n}/b_{1, n}$, $\, a_{2, n}/b_{2, n}$, $\, \cdots, \, a_{2^{n-1}, n}/b_{2^{n-1}, n} = 1/1$, $\, \cdots, \, a_{2^n, n}/b_{2^n, n} = 1/0 \}$, then it is easy to check that 
\begin{itemize}
\item[(a)] $a_{i, n}/b_{i, n} = b_{2^{n}-i, n}/a_{2^{n}-i, n}$ \, for all \, $0 \le i \le 2^{n-1}$; 

\item[(b)] $a_{i, n}/b_{i, n} + a_{2^{n-1}-i, n}/b_{2^{n-1}-i, n} = 1$ \, for all \, $0 \le i \le 2^{n-1}$; 

\item[(c)] $\phi(a_{2^n+i, n}/b_{2^n+i, n}) = a_{i, n}/b_{i, n}$ \, for all \, $0 \le i \le 2^{n-1}$ \, (by (a) and (b)); 

\item[(d)] $\phi(a_{i, n+1}/b_{i, n+1}) =$ (by (a)) $b_{2^n-i, n+1}/a_{2^n-i, n+1} - 1 = a_{2^n-i, n}/b_{2^n-i, n}$ \, for all \, $0 \le i \le 2^n$; 
\end{itemize} 

Now for each positive integer $n$, we define a continuous piecewise linear map $h_n(x)$ from $[0, \infty]$ into $[0, 1]$ by putting $h_n(a_{i, n}/b_{i, n}) = i/2^n$ for all $0 \le i \le 2^n - 1$, $h_n(\infty) = h_n(a_{2^n, n}/b_{2^n, n}) = (2^n - 1)/2^n$, and "connecting the dots".  It is clear that the sequence $< h_n(x) >$ converges uniformly on the non-compact metric space $[0, \infty)$ to a strictly increasing continuous map $h(x)$ with $\lim_{x \to \infty} h(x) = 1$.  Let $h(\infty) = 1$.  Then $h(x)$ is a homeomorphism from $[0, \infty]$ onto $[0, 1]$ and it follows from (c) and (d) above that $h \circ \phi = f \circ h$ on each $F_n$ and so on $[0, \infty]$ by continuity, where $f$ is defined above.  Since $\phi$ is topologically conjugate to $f$ and $f$ has topological entropy $\log \lambda$, where $\lambda$ is the unique positive zero of the polynomial $x^3 - 2x - 1$, so has $\phi$.  We Remark that the map $h$ can also be obtained as the modified Minkowski map $h : [0, \infty] \to [0, 1]$ defined \cite{g} by putting $h(\frac 01) = 0, h(\frac 10) = 1$, and $h(\frac {a+c}{b+d}) = \frac 12[h(\frac ab) + h(\frac cd)]$ whenever $\frac ab$ and $\frac cd$ (both are in lowest terms) are two consecutive fractions in the modified Farey sequences and $h(\frac ab)$ and $h(\frac cd)$ are already defined, and extending $h$ continuously to the whole interval $[0, \infty]$.  

Now since $I_\gamma$ consists of exactly one point for each $\gamma = \gamma_0\gamma_1\gamma_2 \cdots \in \Gamma$ and $\phi(I_{\gamma_0\gamma_1\gamma_2\cdots \gamma_{n_i}}) = I_{\sigma(\gamma_0\gamma_1\gamma_2\cdots \gamma_{n_i})}$, wherever $I_{\gamma_0\gamma_1\gamma_2\cdots \gamma_{n_i}}$ is defined,  it is easy to see that if $\gamma$ and $\gamma(m)$, $m \ge 1$, are points in $\Gamma$ such that $\lim_{m \to \infty} \gamma(m) = \gamma$, then $\lim_{m \to \infty} x_{\gamma(m)} = x_\gamma$.  It is also easy to see that $I_{\ov 0} = \{ (\sqrt{5} - 1)/2 \}$ where $z = (\sqrt{5} - 1)/2$ is the unique fixed point of $\phi(x)$, $I_{1 \ov 0} = \{ (3 + \sqrt{5})/2 \}$, $I_{\overline {010}} = \{ 0 \}$, $I_{\overline {001}} = \{ 1 \}$, and $I_{\overline {100}} = \{ \infty \}$.  This fact will be needed below.  
\bigskip

\noindent
{\bf 5. $\phi(x)$ has bounded uncountable invariant 2-scrambled sets in $[0, 3]$}  

\noindent
The existence of bounded uncountable {\it invariant} $\delta$-scrambled sets for some $\delta > 0$ follow easily from Theorem 3 of \cite{du1}.  Here we use symbolic dynamics to give a different proof.  For each $\beta = \beta_0\beta_1\beta_2 \cdots \in \Sigma_2$, let $\mu_\beta = (\mu_\beta)_0 (\mu_\beta)_1 (\mu_\beta)_2 \cdots$ be a point in $\Gamma$ defined by $$\mu_\beta = 0^{5!} \, A_\beta(5!) \, A_\beta(6!) \, A_\beta(7!) \, \cdots,$$ where $0^2 = 00, 0^3 = 000$, etc., and for each $k \ge 5$, $A_\beta(k!) = (\mu_\beta)_{k!}(\mu_\beta)_{k!+1}(\mu_\beta)_{k!+2} \cdots$ $(\mu_\beta)_{(k+1)!-1}$ is the concatenation of the following $k$ strings of 0's and 1's, each of length $k!$, $$01 0^{k!-2}, \,\, 0\beta_0 0^{k!-2}, \,\, 0\beta_1 0^{k!-2}, \,\, 0\beta_2 0^{k!-2}, \,\, \cdots, \,\, 0\beta_{k-2} 0^{k!-2}.$$ Let $W = \{ \phi^n(x_{\mu_\beta}) : \beta \in \Sigma_2$, $n \ge 0 \}$.  Then it is clear that $W \, (\, \subset I_{000} \cup I_{001} \cup I_{0100} \cup I_{1000} = [0, 1/3] \cup [1/2, 1] \cup [2, 3] \,)$ is a bounded uncountable invariant subset of $[0, 3]$.  If $\eta = \eta_0\eta_1\eta_2 \cdots$ and $\xi = \xi_0\xi_1\xi_2 \cdots$ are any two points (need not be distinct) in $\Sigma_2$ and if $i \ge 1$ and $m \ge 0$ are fixed integers, then, for each $k > i+m+5$, we have 
$$\sigma^{k!+1}(\mu_\eta) = 1 0^{k!-2} \cdots \quad \text{and} \quad \sigma^{k!+2}(\mu_\eta) = 0^{k!-2} \cdots,$$ 
$$\,\,\, \sigma^{k!+1}(\sigma^i(\mu_\xi)) = 0^{k!-1-i} \cdots \quad \text{and} \quad \sigma^{k!+2}(\sigma^i(\mu_\xi)) = 0^{k!-2-i} \cdots,$$ and 
$$\,\,\, \sigma^{(m+2)k!+1}(\mu_\eta) = \eta_m 0^{k!-2} \cdots \quad \text{and} \quad \sigma^{(m+2)k!+1}(\mu_\xi) = \xi_m 0^{k!-2} \cdots.$$   

As $k$ tends to $\infty$, the sequence $<1 0^{k!-2} \cdots>$ tends to the point $1 \bar 0$ and the sequence $<0^{k!-2-i} \cdots>$ tends to the point $\bar 0$.  Since $\lim_{n \to \infty} I_{0^{k!-2-i} \cdots} = I_{\ov 0} = \{ (\sqrt{5} - 1)/2 \} = \{ x_{\ov 0} \} = \{$the unique fixed point of $\phi \}$ and $\lim_{n \to \infty} I_{1 0^{k!-2} \cdots} = I_{1 \ov 0} = \{ (3 + \sqrt{5})/2 \} = \{ x_{1 \ov 0} \} = \{$the unique inverse image in $(1, \infty)$ of the unique fixed point of $\phi \}$, we obtain $$\limsup_{n \to \infty} |\phi^n(x_{\mu_\eta}) - \phi^n(\phi^i(x_{\mu_\xi}))| \ge |x_{1 \ov 0} - x_{\ov 0}| = 2 \quad \text{and} \quad \liminf_{n \to \infty} |\phi^n(x_{\mu_\eta}) - \phi^n(\phi^i(x_{\mu_\xi}))| = 0.$$  Similarly, if $\eta_m \ne \xi_m$, then we have 
$$\limsup_{n \to \infty} |\phi^n(x_{\mu_\eta}) - \phi^n(x_{\mu_\xi})| \ge 2 \quad \text{and} \quad 
\liminf_{n \to \infty} |\phi^n(x_{\mu_\eta}) - \phi^n(x_{\mu_\xi})| = 0.$$  Therefore, we have proved the following result:  

\noindent
{\bf Theorem 1.} {\it Let $\phi(x) = |1 - 1/x|$ for all $x > 0$.  Then $W \, (\, \subset [0, 3]\,)$ is a bounded uncountable invariant 2-scrambled set for $\phi(x)$.}
\bigskip

\noindent
{\bf 6. $\phi(x)$ has dense unbounded uncountable invariant $\infty$-scrambled sets}

\noindent
Let $k \ge 5$ be a fixed integer.  We call any finite sequence $\beta_0\beta_1\beta_2 \cdots \beta_k$ of 0's and 1's admissible if $\beta_i = 1$ and $i < k$ then $\beta_{i+1} = 0$.  There are countably infinitely many such admissible finite sequences of 0's and 1's.  Let $B_1, B_2, \cdots, B_n, \cdots$ be an enumeration of these admissible finite sequences with $B_n = \beta_{n,0}\beta_{n,1}\beta_{n,2} \cdots \beta_{n, k_n-1}$, $n \ge 1$.  Let $m_1, m_2, \cdots$ be a strictly increasing sequence of positive integers such that $(m_i)! > (m_{i-1})! + k_{i-1} + 1$ for all $i \ge 2$.  Let $\alpha =$ $\alpha_0\alpha_1\alpha_2 \cdots$ be a point in $\Sigma_2$ defined by putting $\alpha_{(m_i)!} \alpha_{(m_i)! + 1} \alpha_{(m_i)! + 2} \cdots \alpha_{(m_i)! + k_i-1} = B_i$ for all $i \ge 1$ and $\alpha_\ell = 0$ elsewhere.  Then it is clear that $\alpha$ is a {\it totally transitive point} of $\sigma$ in $\Gamma = \{ \gamma = \gamma_0\gamma_1\gamma_2 \cdots \in \Sigma_2 : $ if $\gamma_i = 1$ then $\gamma_{i+1} = 0 \}$, i.e., for each integer $n \ge 1$, the orbit of the point $\alpha$ with respect to the map $\sigma^n$ is dense in $\Gamma$.  

Let $\alpha = \alpha_0\alpha_1\alpha_2 \cdots$ be a {\it totally transitive point} of $\sigma$ in $\Gamma$ (and so if $I_\alpha = \{ x_\alpha \}$ then $x_\alpha$ is a {\it totally transitive point} of $\phi$ in $(0, \infty)$).  Let $X = \{ x(1), x(2), x(3), \cdots \}$ be any countably infinite subset of {\it irrational} points in $(0, \infty)$.  For each integer $m \ge 1$, we can write $\{ x(m) \} = I_{\gamma^{(m)}}$ for an unique $\gamma^{(m)} = \gamma^{(m)}_0\gamma^{(m)}_1\gamma^{(m)}_2 \cdots$ in $\Gamma$.  For any integers $m \ge 1$ and $5 \le i < j$ such that $j-i+1$ is a multiple of 3, let 
$$C(x(m), \, i : j) = \gamma^{(m)}_i\gamma^{(m)}_{i+1} \cdots \gamma^{(m)}_{j-1}\,\, 0 \quad \text{and} \,\,\,\qquad\qquad\qquad  \qquad\qquad\qquad\qquad\qquad\qquad$$ 
$$
C^*(x(m), \, i : j) = \xi^{(m)}_i\xi^{(m)}_{i+1}\xi^{(m)}_{i+2} \cdots \xi^{(m)}_j = \begin{cases}
      0\, 100 \, 100 \, 100 \, \cdots 100 \, 10, \quad \text{if \, $\gamma^{(m)}_i = 1$}, \\
      100 \, 100 \, 100 \, \cdots 100, \qquad\quad \text{otherwise}. \\
      \end{cases}
$$
Note that, in $C(x(m), \, i : j)$ and in $C^*(x(m), \, i : j)$, the last element is 0 so when we concatenate $C$'s with $C$'s or with $C^*$'s, we won't have the finite string $\cdots 11 \cdots$.  For simplicity, for $b = 0$ or 1, we write $(b 00)^1 = b 00, (b 00)^2 = b 00 \, b 00$, and so on.  

For any $\beta = \beta_0 \beta_1 \beta_2 \cdots$ in $\Sigma_2$, define a new point $\tau_\beta = \tau_\beta(\alpha,X)$ in $\Gamma$ as follows and let $$Y = \{ \phi^n(x_{\tau_\beta}) : \beta \in \Sigma_2, \, n \ge 0 \},$$ 
where $\tau_\beta = (\tau_\beta)_0(\tau_\beta)_1(\tau_\beta)_2 \cdots = \alpha_0\alpha_1\alpha_2 \cdots \alpha_{5!-2}0$ $A_\beta(5!)$ $A_\beta(6!)$ $A_\beta(7!)$ $\cdots$,
and, for each $k \ge 5$, $A_\beta(k!) = (\tau_\beta)_{k!}(\tau_\beta)_{k!+1}(\tau_\beta)_{k!+2} \cdots (\tau_\beta)_{(k+1)!-1}$ is the concatenation of the following $k$ strings of 0's and 1's, each of length $k!$,

\qquad $\alpha_0\alpha_1\alpha_2 \cdots \alpha_{k!-1}$ 

\qquad $0^{k!/4} (100)^{k!/12} (001)^{k!/12} (010)^{k!/12}$ 

\qquad $(\beta_0 0 0)^{(k-1)!/3}(\beta_1 0 0)^{(k-1)!/3} \cdots (\beta_{k-1} 0 0)^{(k-1)!/3}$ 

\qquad $B(x(1), 4k!)$ \qquad $B(x(2), 5k!)$ \qquad $B(x(3), 6k!)$ \qquad $\cdots$ \qquad $B(x(k-3), k \cdot k!)$,

\noindent
where, for $1 \le i \le k-3$, $B(x(i), (3+i)k!)$ is the concatenation of the following $2k$ strings of 0's and 1's, each of length $\frac 12 (k-1)!$, 

\qquad $C(x(i), (3+i)k!: (3+i)k!+[{\f 12}(k-1)!-1])$

\qquad \qquad\,\, $\cdots$ 

\qquad $C(x(i), (3+i)k!+(j-1)[{\f 12}(k-1)!-1]: (3+i)k!+j[{\f 12}(k-1)!-1])$

\qquad \qquad\,\, $\cdots$ \,\, \quad 

\qquad $C(x(i), (3+i)k!+(k-1)[{\f 12}(k-1)!-1]: (3+i)k!+k[{\f 12}(k-1)!-1])$

\qquad $C^*(x(i), (3+i)k!+{\f 12}k!: (3+i)k!+{\f 12}k!+[{\f 12}(k-1)!-1)])$ 

\qquad \qquad\,\, $\cdots$ 

\qquad $C^*(x(i), (3+i)k!+{\f 12}k!+(j-1)[{\f 12}(k-1)!-1]: (3+i)k!+{\f 12}k!+j[{\f 12}(k-1)!-1])$

\qquad \qquad\,\, $\cdots$ \,\, \quad 

\qquad $C^*(x(i), (3+i)k!+{\f 12}k!+(k-1)[{\f 12}(k-1)!-1]: (3+i)k!+{\f 12}k!+k[{\f 12}(k-1)!-1])$.

In the following, we shall use the fact that $I_{\overline {100}} = \{ \infty \}$, $I_{\overline {001}} = \{ 1 \}$, $I_{\overline {010}} = \{ 0 \}$, and if $\gamma$ and $\gamma^{(n)}$, $n \ge 1$, are points in $\Gamma$ such that $\lim_{n \to \infty} \gamma^{(n)} = \gamma$ then $\lim_{n \to \infty} \phi(x_{\gamma^{(n)}}) = x_\gamma$, where $I_{\gamma^{(n)}} = \{ x_{\gamma^{(n)}} \}$ and $I_\gamma = \{ x_\gamma \}$.  Now in the expansion of $\tau_\beta$, $\beta \in \Sigma_2$, 
\begin{itemize}
\item[(1)] there are infinitely many strings $\alpha_0\alpha_1\alpha_2 \cdots \alpha_{k!-1}$, $k \ge m$, which imply, for each \,$i \ge 1$, the denseness of the orbit  (under $\phi^i$) of the point $x_{\tau_\beta}$ in $(0, \infty)$, 

\item[(2)]
there are infinitely many strings $(\beta_0 0 0)^{(k-1)!/3}(\beta_1 0 0)^{(k-1)!/3} \cdots (\beta_{k-1} 0 0)^{(k-1)!/3}, \,\, k \ge 5$ which imply that $$\limsup\limits_{n \to \infty} |\phi^n(x_{\tau_\beta}) - \phi^n(x_{\tau_\eta})| = \infty \,\, \text{for} \,\,\, \beta \ne \eta \,\,\, \text{in} \,\,\, \Sigma_2,$$

\item[(3)]
for any $i \ge 1$, there are infinitely many strings $0^{k!/4} (100)^{k!/12} (001)^{k!/12} (010)^{k!/12}$, $k \ge i+5$,  containing $$0^i (100)^{k!/12} (001)^{k!/12} (010)^{k!/12}$$ which imply that, for any integer $i \ge 1$, $$\limsup\limits_{n \to \infty} |\phi^n(x_{\tau_{\beta}}) - \phi^n(\phi^i(x_{\tau_{\eta}}))| = \infty \,\,\, \text{for any} \,\,\, \beta \,\,\, \text{and} \,\,\, \eta \,\,\, \text{in} \,\,\, \Sigma_2,$$ and since this string also contains long strings of $0$'s, i.e., $0^{k!/4}$, $k \ge 5$, we obtain $$\liminf\limits_{n \to \infty} |\phi^n(x_{\tau_{\beta}}) - \phi^n(\phi^i(x_{\tau_{\eta}}))| = 0 \,\,\, \text{for any} \,\,\, \beta, \eta \,\,\, \text{in} \,\,\, \Sigma_2 \,\,\, \text{and} \,\,\, i \ge 0,$$

\item[(4)]
for any positive integers $i$ and $j$, there are infinitely many strings $$C(x(i), (3+i)k!+(j-1)[{\f 12}(k-1)!-1]: (3+i)k!+j[{\f 12}(k-1)!-1]),\,\, k > j,$$ which imply that $$\liminf\limits_{n \to \infty} |\phi^n(x(i)) - \phi^n(\phi^{j-1}(x_{\tau_\beta}))| = 0 \,\,\, \text{for all} \,\,\, \beta \,\,\, \text{in} \,\,\, \Sigma_2 \,\,\, \text{and} \,\,\, j \ge 1,$$

\item[(5)]
for any positive integers $i$ and $j$, there are infinitely many strings $$C^*(x(i), (3+i)k!+{\f 12}k!+(j-1)[{\f 12}(k-1)!-1]: (3+i)k!+{\f 12}k!+j[{\f 12}(k-1)!-1]), \,\, k > j,$$ which, since $x_{0\gamma_1\gamma_2 \cdots} \le 1 \le x_{1\gamma_1 \gamma_2 \cdots}$ and $\lim_{k \to \infty} I_{(100)^k} = \{ x_{\ov {100}} = \infty \}$, imply that $$\limsup\limits_{n \to \infty} |\phi^n(x(i)) - \phi^n(\phi^{j-1}(x_{\tau_\eta}))| = \infty \,\,\, \text{for all} \,\,\, \beta \,\,\, \text{in} \,\,\, \Sigma_2 \,\,\, \text{and} \,\,\, j \ge 1.$$
\end{itemize}

Now, if $r$ is a rational number in $[0, \infty)$, then, by \cite{se1}, for some integer $m \ge 0$, $\phi^m(r) = 0$ which is a period 3 point of $\phi$.  For each $\beta$ in $\Sigma_2$, the iterates of the point $x_{\tau_\beta}$ approach the fixed point $z = (\sqrt{5} - 1)/2$ infinitely often and stay close to it for a while each time.  This implies that $\limsup_{n \to \infty} |\phi^n(r) - \phi^n(x_{\tau_\beta})| = \infty$.  Furthermore, since, for all $k > 12$, $$\sigma^{3(k!/3+k!/12)}(\tau_\beta) = \sigma^{k!+k!/4}(\tau_\beta) = (100)^{k!/12}(001)^{k!/12}(010)^{k!/12} \cdots,$$ $$\sigma^{3(k!/3+k!/12+k!/36)}(\tau_\beta) = \sigma^{k!+k!/4+k!/12}(\tau_\beta) = (001)^{k!/12}(010)^{k!/12} \cdots,$$ and $$\sigma^{3(k!/3+k!/12+k!/36+k!/36)}(\tau_\beta) = \sigma^{k!+k!/4+k!/12+k!/12}(\tau_\beta) = (010)^{k!/12} \cdots,$$ the sequence $< \sigma^{3n}(\tau_\beta) >$ approaches each of the three points $\overline {100}$, $\overline {001}$ and $\overline {010}$ infinitely often and so the sequence $< \phi^{3n}(x_{\tau_\beta}) >$ approaches each of the three points $\infty, 1, 0$ infinitely often.  Since, for some $m \ge 1$, $\phi^m(r) = 0$ and $0$ is a period-3 point of $\phi(x)$, this implies that $\liminf_{n \to \infty} |\phi^n(r) - \phi^n(x_{\tau_\beta})| = 0$ (if $f^{3n}(r) = \infty$ and $\phi^{3n}(x_{\tau_\beta}) \approx \infty$ then $\phi^{3n+1}(r) = 1$ and $\phi^{3n+1}(x_{\tau_\beta}) \approx 1$).  Therefore, we have shown the following result:  

\noindent
{\bf Theorem 2.}
{\it Let $\phi(x) = |1 - 1/x|$ for all $x > 0$.  Then for any given countably infinite subset $X$ of points (rational or irrational) in $(0, \infty)$, there exists a dense unbounded uncountable invariant $\infty$-scrambled set $Y$ of totally transitive irrational points in $(0, \infty)$ such that, for any $x$ in $X$ and any $y$ in $Y$, we have $$\limsup_{n \ra \infty} |\phi^n(x) - \phi^n(y)| = \infty \quad \text{and} \quad \liminf_{n \ra \infty} |\phi^n(x) - \phi^n(y)| = 0.$$}
\indent
In the above result, we consider $\phi(x)$ as a continuous map from the compact topological (but not metric) space $[0, \infty]$ onto itself under the convention that $|\frac 10| = \infty$, $\frac 1{\pm \infty} = 0$, and $|\pm \infty \, \pm$ any real number$| = \infty$.  It is easy to see that the dynamics of $\phi(x)$ on the rational points in $[0, \infty)$ are not interesting \cite{se1, se2} because they all go to the period-3 orbit $\{ 1, 0, \infty \}$ after a finite number of iterations.  All interesting dynamics occur in the invariant set $\mathbb{R_+} \setminus \mathbb{Q_+}$ of irrational points of $(0, \infty)$.  Consequently, we can consider the dynamics of $\phi(x)$ on the set $\mathbb{R_+} \setminus \mathbb{Q_+}$ which is a metric space under the usual distance metric.  Within this context, Theorem 2 has a very important consequence: given any {\it irrational} point $x$ in $\mathbb{R_+} \setminus \mathbb{Q_+}$, then at just about everywhere in $\mathbb{R_+} \setminus \mathbb{Q_+}$, whether it is close to $x$ or far away from it we can always find an {\it irrational} point $y$ (in the dense set $Y$) whose iterates satisfy $\limsup_{n \ra \infty} |\phi^n(x) - \phi^n(y)| = \infty$ and $\liminf_{n \ra \infty}$ $|\phi^n(x) - \phi^n(y)| = 0$.  This demonstrates the true nature of chaos \cite{du2, v1, v2}, i.e., not only nearby points will separate (sensitivity) and converge infinitely often but even {\it far apart} points will also converge and separate infinitely often.    
\bigskip

\noindent
{\bf Acknowledgement}

\noindent
This work was partially supported by the National Science Council NSC 97-2119-M-001-004.

\end{document}